\documentclass{elsart}
\usepackage{amsmath,amsfonts,amssymb,indentfirst}

\newcommand{\fb}{\mathbf}
\newcommand{\<}{\kern.0833em}
\def\r#1{{\rm #1}}
\newcommand{\coLm}{\hspace{.05em}{\varinjlim\hspace{.1em}}}
\newcommand{\Lm}{\hspace{.05em}{\varprojlim\hspace{.05em}}}
\newcommand{\ObE}{\mathrm{Ob}(\fb E)}
\newcommand{\D}{{\rm down}}

\makeatletter
\newcommand{\xlabel}{\stepcounter{equation}
  \gdef\@currentlabel{\p@equation\theequation}{\rm(\@currentlabel)}}
\makeatother
\newenvironment{xlist}
  {\begin{list}{\xlabel}
    {\setlength{\rightmargin}{20pt}
     \setlength{\leftmargin}{37pt}
     \setlength{\labelsep}{20pt}
     \setlength{\labelwidth}{20pt}}}
  {\end{list}}

\begin{document}
\begin{frontmatter}
\title{Direct limits and fixed point sets%
\thanksref{foot}}
\thanks[foot]{
Preprint versions:
\ http://math.berkeley.edu/\protect\linebreak[0]%
{$\!\sim$}gbergman/papers/%
dirlimfix.\{tex,dvi,ps\},~
~arXiv:math.CT/0306127.}
\author{George M. Bergman}
\address{University of California,
Berkeley, CA 94720-3840,
USA\\
{\rm gbergman@math.berkeley.edu}
}

\begin{abstract}
For which groups $G$ is it true that whenever one forms a direct
limit of left $G$-sets, $\coLm_{i\in I}\,X_i,$ the set of its fixed
points, $(\coLm_I\,X_i)^G,$ can be obtained as the direct limit
$\coLm_I(X_i^G)$ of the fixed point sets of the given $G$-sets?
An easy argument shows that this is the case if and only if $G$ is
finitely generated.

If we replace ``group $G$'' by ``monoid $M$'', the answer is the
less familiar condition that the improper left congruence on $M$
be finitely generated; equivalently, that $M$ be finitely generated
under multiplication and ``right division''.

Replacing our group or monoid with a small category $\fb E,$
the concept of a set on which $G$ or $M$ acts with that
of a functor $\fb E\rightarrow \fb{Set},$ and the
fixed point set of an action with the limit of a functor,
a criterion of a similar nature is proved.
Specialized criteria are obtained in the cases where $\fb E$
has only finitely many objects and where $\fb E$ is a
(generally infinite) partially ordered set.

If one allows the codomain category $\fb{Set}$ to be replaced with other
categories, and/or allows direct limits to be replaced with other
classes of colimits, one enters a vast area open
to further investigation.
\end{abstract}
\end{frontmatter}

\section{Introduction.}\label{outline}
Although the next three sections, concerning fixed point sets of
group and monoid actions, require no familiarity
with category theory, I will (with apologies to the non-categorical
reader) frame this introduction in category-theoretic terms.

It is a familiar observation that ``left universal constructions
respect left universal constructions and right universal constructions
respect right universal constructions'' \cite[\S\S7.7-7.8]{245}.
Thus, when one takes a limit of limits, or a colimit of
colimits (in a context where the relevant limits or colimits
all exist), one can reverse the order of the two limit operations,
or of the two colimit operations, without changing the result.
In contrast, left and right universal constructions do not in general
respect one another.
(For instance, the free group on a direct product set $X\,{\times}\,Y$
is not isomorphic to the direct product of the free group on $X$
and the free group on~$Y).$

But there are classes of cases where, anomalously,
certain limits commute with certain colimits.
For instance, given directed systems of sets $(X_i)_I$ and
$(Y_i)_I$ indexed by the same partially ordered set $I,$
one finds that $\coLm (X_i\,{\times}\,Y_i)\cong (\coLm X_i)\times
(\coLm Y_i).$
Indeed, the fact that we can construct a direct limit of algebras
by putting an algebra structure on the direct limit of their
underlying sets is a consequence of this fact, given that
algebra operations on $X$ are set maps $X{\times}\cdots{\times}X\to X.$

This note investigates the question of which small categories
$\fb E$ have the property that limits of functors from $\fb E$ to
$\fb{Set}$ always commute with direct limits, that is, with
colimits over directed partially ordered sets.
It has been observed
\cite[Thm.~IX.2.1, p.211]{CW}, \cite[Thm.~4.73, p.72]{MK} that
this happens if $\fb E$ is a finite category, i.e., has only
finitely many objects and finitely many morphisms.
More generally, it occurs whenever $\fb E$ has finitely many
objects and finitely {\em generated} morphism-set
(\cite[Prop.~7.9.3]{245} $=$ Corollary~\ref{CfgE} below).
The result of \S\ref{group} (first paragraph of
the above abstract) is equivalent to the statement that
if $\fb E$ is a one-object category whose morphisms form a {\em group},
this finite generation condition is necessary as well as sufficient.

In a general one-object category $\fb E,$ the morphisms form a
monoid $M.$
By the result noted above, finite generation of $M$ is
sufficient
for the construction of limits over $\fb E$ (i.e., fixed-point sets of
$\!M\!$-sets) to commute with that of direct limits, but in this case it
is not necessary.
In \S\ref{moncrit} we obtain two
criteria each of which is necessary as well as sufficient.
We find in \S\ref{cat} that one of these, finite generation
of the improper left congruence on $M,$ when reformulated as
finite presentability of the trivial $\!M\!$-set, generalizes to
arbitrary small categories $\fb E,$ while the other, finite
generation of $M$ under multiplication and ``right division'',
generalizes nicely to categories $\fb E$ with finitely many objects.

In \S\ref{poset} we examine the case where $\fb E$ is the category
$J_\fb{cat}$ induced by a partially ordered set $J,$ and
translate our general criterion into a condition on~$J.$
Half of the condition we get can be stated in familiar language: It says
that the set of minimal elements of $J$ is finite, and every
element lies above a minimal element.
(This is in fact necessary and sufficient for
the comparison maps associated with our limits and colimits to be
{\em injective} in all cases;
it is also {\em necessary} for them always to be surjective.)
The remaining condition appears to be new.
In language which we shall define, it says that the set of elements
of $J$ ``critical'' with respect to the minimal elements is
finite, and that these critical elements
``gather'' all minimal elements under every element of~$J.$

Note that the results of this paper only concern
limits and colimits of functors to $\fb{Set};$
the behavior of functors to other
categories can be strikingly different.
For instance \cite[Exercise~7.9.5]{245},
in $\fb{Set}^{\rm op},$ direct limits do not in general commute with
equalizers, though equalizers are limits over a certain finite category;
but they do commute with not necessarily finite small products;
so we have both negative and positive deviations from the behavior of
$\!\fb{Set}\!$-valued functors.
Clearly, it would be interesting to investigate more classes of cases
of commutativity between limits and colimits: for functors with
codomains other than $\fb{Set},$ and for
colimits over categories other than directed partially ordered sets.
If we fix one of the three variables -- the small category over which
we take limits, the small category over which we take colimits,
and the codomain category -- then we get a Galois connection
\mbox{\cite[\S5.5]{245}} on the other two, and can study the
resulting closure operators.
The exercises in \cite[\S{IX.2}]{CW} and the results and
exercises at the end of \cite[\S7.9]{245} give scattered results
along these lines, but for the most part,
the topic seems wide open for study!

I am indebted to Birge Huisgen-Zimmermann and Ken Goodearl for
organizing the gathering at which I first spoke about some of these
results, to Max Kelly, Arthur Ogus, and Boris Schein for references to
the literature, and to the referee for several helpful corrections
and suggestions.

The present note has various possible audiences, ranging from any
mathematician who uses direct limits, to the specialist in
semigroups or categories or partially ordered sets.
I hope the reader will be patient with my reviewing
details that may be familiar to him or her, and also with my following,
in \S\ref{moninit}, a somewhat leisurely path of motivation
to the results on monoids.

\section{Direct limits and group actions.}\label{group}

Recall that a partially ordered set $(I,\leq)$ is said
to be {\em directed} if for every pair of elements $i,j\in I,$
there exists $k\in I$ majorizing both,
i.e., satisfying $k\geq i$ and $k\geq j.$
A {\em directed system of sets} means a family of sets
$(X_i)_{i\in I}$ indexed by a nonempty directed partially
ordered set $I,$ and given with connecting maps
$\alpha_{i,j}\<{:}\ X_i\to X_j$ $(i\leq j)$ such that
each $\alpha_{i,i}$ is the identity map of $X_i,$ and whenever
$i\leq j\leq k,$ one has $\alpha_{i,k}=\alpha_{j,k}\,\alpha_{i,j}.$
(So a more complete notation for the directed
system is $(X_i,\alpha_{i,j})_{i,j\in I}.)$

In this situation one has the concept of the {\em direct limit} of
the given system.
This is constructed by forming the disjoint union $\bigsqcup_IX_i,$
and dividing out by the least equivalence relation $\sim$
such that $x\sim \alpha_{i,j}(x)$
whenever $x\in X_i$ and $i\leq j.$
Denoting the resulting set
$\coLm_I\,X_i,$ and writing $[x]$ for the equivalence class
therein of $x\in\bigsqcup_IX_i,$ we get, for each $j\in I,$ a map
$\alpha_{j,\infty}\<{:}\ X_j\to \coLm_I\,X_i$
taking $x\in X_j$ to $[x].$
The characterization of $\coLm_I\,X_i$ that we will use here
is that it is a set given with maps
$\alpha_{j,\infty}\<{:}\ X_j\to \coLm_I\,X_i$ for
each $j\in I,$ such that every
element of $\coLm_I\,X_i$ is of the form $\alpha_{j,\infty}(x)$
for some $j\in I,$ $x\in X_j,$ and such that

\begin{xlist}\item\label{eqiff}
$\alpha_{i,\infty}(x)\,=\,\alpha_{j,\infty}(y)\,$ if and only if there
exists $k\geq i,\,j$ such that $\,\alpha_{i,k}(x)\,=\,\alpha_{j,k}(y).$
\end{xlist}

Property~(\ref{eqiff}) is easily deduced from the above
construction of $\coLm_I\,X_i,$ using the directedness of $I.$
Note that it includes the relations

\begin{xlist}\item\label{ijinf}
$\alpha_{i,\infty}(x)\,=\,\alpha_{j,\infty}(\alpha_{i,j}(x))\quad
(i\leq j\in I,~x\in X_i)$
\end{xlist}
corresponding to the generators of the equivalence relation
in that construction.

If $G$ is a group, then a directed system of left $\!G\!$-sets
means a directed system $(X_i,\alpha_{i,j})_{i,j\in I}$ of sets, such
that each $X_i$ is given with a left action of $G,$ and each of the
connecting maps $\alpha_{i,j}$ is a morphism of $\!G\!$-sets (a
$\!G\!$-equivariant map).
Henceforth we will generally omit the qualifier ``left''.
Given such a directed system, it is easy to verify
that $\coLm_I\,X_i$ admits a unique $\!G\!$-action making
the maps $\alpha_{i,\infty}$ morphisms of $\!G\!$-sets, i.e., such that

\begin{xlist}\item\label{g&ioo}
$g\,\alpha_{i,\infty}(x)\,=\,\alpha_{i,\infty}(gx)\quad
(g\in G,\ i\in I,\ x\in X_i).$
\end{xlist}

For any $\!G\!$-set $X,$ let us write
\begin{xlist}\item[]
$X^G\,=\,\{x\in X~|~(\<\forall g\in G)\,\ gx=x\}$
\end{xlist}
for the fixed-point set of the action.
If $(X_i,\alpha_{i,j})_{i,j\in I}$ is a directed
system of $\!G\!$-sets, we see that each map $\alpha_{i,j}$ carries the
fixed set $X^G_i$ into $X^G_j.$
Writing $\beta_{i,j}$ for the restriction of $\alpha_{i,j}$ to a
map $X^G_i\to X^G_j,$ we thus get a directed system of sets
$(X^G_i,\beta_{i,j}),$ and we can form its direct
limit $\coLm_I\,X^G_i.$

It is now straightforward to verify that one has a map

\begin{xlist}\item\label{iotaG}
$\iota:\ \coLm_I\,X^G_i\longrightarrow\,(\coLm_I\,X_i)^G,$ \ \ defined
by \ \ %
$\iota(\beta_{i,\infty}(x))\,=\,\alpha_{i,\infty}(x)$ $(x\in X^G_i).$
\end{xlist}

\begin{thm}\label{TGfixed}
If $G$ is a group, $I$ a directed partially ordered
set, and\linebreak[4] $(X_i,\alpha_{i,j})_{i,j\in I}$ a directed
system of $\!G\!$-sets, then the set-map $\iota$ of~{\rm(\ref{iotaG})}
is one-to-one.

Moreover, for any group $G,$ the following conditions are equivalent:
\begin{xlist}
\item\label{bijG}
For every directed partially ordered set $I$ and
directed system\linebreak[4] $(X_i,\alpha_{i,j})_{i,j\in I}$ of
$\!G\!$-sets, the set-map $\iota$ of~{\rm(\ref{iotaG})} is bijective.
\end{xlist}
\begin{xlist}
\item\label{Gfg}
$G$ is finitely generated.
\end{xlist}
\end{thm}\begin{pf}
The assertion of the first
sentence follows from~(\ref{eqiff}) and the fact
that the maps $\beta_{i,j}$ are restrictions of the $\alpha_{i,j}.$

To see that~(\ref{Gfg}) implies~(\ref{bijG}),
let $\{g_1,\dots,g_n\}$ be a finite generating set for $G,$ and
consider any element of $(\coLm_I\,X_i)^G,$ which we may write
$\alpha_{i,\infty}(x)$ for some $i\in I$ and $x\in X_i.$
The element $x\in X_i$ may not itself be fixed
under $G,$ but by assumption, for every $g\in G$ we have
$g\,\alpha_{i,\infty}(x)=\alpha_{i,\infty}(x),$ in other words,
$\alpha_{i,\infty}(gx)=\alpha_{i,\infty}(x).$
By~(\ref{eqiff}) this means that for each $g\in G$ there
exists $k(g)\geq i$ in $I$ such that
$\alpha_{i,k(g)}(gx)=\alpha_{i,k(g)}(x).$

Since $I$ is directed, we can find a common upper bound
$k$ for $k(g_1),\dots,k(g_n),$ and we see from the
$\!G\!$-equivariance of the maps $\alpha_{k(g_j),k}$ that
$\alpha_{i,k}(x)$ will be invariant under all
of $\{g_1,\dots,g_n\},$ hence will belong to $X^G_k.$
The element $\beta_{k,\infty}(\alpha_{i,k}(x))$ is thus an element
of $\coLm_I\,X^G_i,$ and~(\ref{ijinf}) shows that it
is mapped by~(\ref{iotaG}) to the given
element $\alpha_{i,\infty}(x)\in (\coLm_I\,X_i)^G,$ as required.

Conversely, if $G$ is a non-finitely-generated group, let
$I$ be the set of finitely generated subgroups of $G,$ partially
ordered by inclusion; this is clearly a directed partially ordered set.
For each $H\in I,$ let $X_H$ be the transitive $\!G\!$-set $G/H,$
and define connecting maps by $\alpha_{H_1,H_2}(gH_1) = gH_2$
for $H_1\leq H_2;$ this gives a directed system.
Since each $H\in I$ is a proper subgroup of $G,$ each of the
$\!G\!$-sets $X_H$ satisfies $(X_H)^G = \varnothing,$
so $\coLm_I\,(X_H)^G=\varnothing.$
On the other hand, any two elements $g_1H_1\in X_{H_1}$
and $g_2H_2\in X_{H_2}$ have the same image in
$X_{H_3}$ for any $H_3$ containing $H_1,\; H_2,$ and
$g_1^{-1}g_2\,,$ so $\coLm_I\,X_H$ is the one-point $\!G\!$-set.
Thus $(\coLm_I\,X_H)^G\neq\varnothing,$ and~(\ref{bijG}) fails.\qed
\end{pf}

{\em Digression.}
One may ask whether~(\ref{bijG}) is equivalent to the corresponding
statement with $I$ restricted to be the
set $\Nset$ of natural numbers with the usual ordering~${\leq\<},$
this being the kind of direct limit one generally first learns about.
If we call this weakened condition~(\ref{bijG}$\!_{\Nset}\!$),
I claim the proof of Theorem~\ref{TGfixed} may be adapted to
show that~(\ref{bijG}$\!_{\Nset}\!$) is equivalent to
\begin{xlist}
\item[\ \ (\ref{Gfg}$\!_{\Nset}\!$)]
Every chain $H_0\leq H_1\leq\dots$ of subgroups of $G$
indexed by ${\Nset}$ and having union $G$ is eventually constant.
\end{xlist}
Indeed, suppose $G$ is a group for which~(\ref{bijG}$\!_{\Nset}\!$)
fails, so that we have a directed system $(X_i)_{i\in\Nset}$ and
an element $\alpha_{j,\infty}(x)\in(\coLm_{\Nset}\,X_i)^G$
which is not in the image of $\iota.$
Then no $\alpha_{j,k}(x)$ lies in $X^G_k,$ and letting
$H_i$ be the isotropy subgroup of $\alpha_{j,j+i}(x)$ for each $i,$
it is easy to see that these subgroups give
a counterexample to~(\ref{Gfg}$\!_{\Nset}\!$).
Conversely, if we have a counterexample
to~(\ref{Gfg}$\!_{\Nset}\!$), then setting $X_i=G/H_i$
gives a counterexample to~(\ref{bijG}$\!_{\Nset}\!$).

But are there any groups that satisfy~(\ref{Gfg}$\!_{\Nset}\!$)
and not~(\ref{Gfg})?
Clearly~(\ref{Gfg}$\!_{\Nset}\!$)
cannot hold in any countable non-finitely-generated group.
It will also fail in any group which admits a homomorphism onto
a group in which it fails, from which one can show that it fails
in any non-finitely-generated abelian
group \cite[paragraph following Question~8]{Sym_Omega:1}.
However, examples are known of uncountable nonabelian groups
that satisfy~(\ref{Gfg}$\!_{\Nset}\!$):
Infinite direct powers of nonabelian simple
groups \cite{SK+JT}, full permutation groups on
infinite sets \cite{MN,Sym_Omega:1}, and others
\cite{MD+RG,MD+WCH,Sh,STcof}.

(Groups satisfying~(\ref{Gfg}$\!_{\Nset}\!$)
but not~(\ref{Gfg}) are said to be of ``uncountable cofinality''.
The same condition on modules has been studied under a surprising
variety of names \cite[p.895, top paragraph]{EKN}.)

\section{Monoid actions -- initial observations.}\label{moninit}

If we replace the group $G$ of the preceding section with a general
monoid $M,$ a large part of the discussion goes over unchanged.
Given a directed system $(X_i,\alpha_{i,j})_{i,j\in I}$ of left
$\!M\!$-sets, we get an $\!M\!$-set structure on $\coLm_I\,X_i,$
and there is a natural map
\begin{xlist}
\item\label{iotaM}
$\iota:\ \coLm_I\,X^M_i\longrightarrow\,(\coLm_I\,X_i)^M$\quad given
by $\iota(\beta_{i,\infty}(x))=\alpha_{i,\infty}(x)$ $(x\in X^M_i),$
\end{xlist}
which is always one-to-one; and again we may ask for which $M$ it
is true that
\begin{xlist}\item\label{bijM}
For every directed partially ordered set $I$ and
directed system $(X_i,\alpha_{i,j})_{i,j\in I}$ of
$\!M\!$-sets, the set-map $\iota$ of~(\ref{iotaM}) is bijective.
\end{xlist}
The argument used in the proof of Theorem~\ref{TGfixed},
(\ref{Gfg})\!$\implies$\!(\ref{bijG}),
shows that a {\em sufficient} condition is
\begin{xlist}\item\label{Mfg}
$M$ is finitely generated.
\end{xlist}
Attempting to reproduce the converse argument, we can say, as before,
that if $M$ is not finitely generated its finitely generated
submonoids $N$ form a directed partially ordered
set; however, there is no concept of
factor-$\!M\!$-set $M/N,$ as would be needed to continue the argument.

And in fact, there exist non-finitely-generated monoids for
which~(\ref{bijM}) holds.
For instance, let $M$ be the multiplicative monoid of any
field $F;$ note that $0\in M.$
Given an element $\alpha_{j,\infty}(x)\in (\coLm_I\,X_i)^M,$
we have $\alpha_{j,\infty}(x)= 0\,\alpha_{j,\infty}(x)=
\alpha_{j,\infty}(0x),$ hence there exists
$k\in I$ such that $\alpha_{j,k}(x)= \alpha_{j,k}(0x).$
We now observe that for every $u\in M$ we have
\begin{xlist}\item[]
$u\,\alpha_{j,k}(x)\ =\ u\,\alpha_{j,k}(0x)\ =\ %
\alpha_{j,k}((u\,0)x)\ =\ \alpha_{j,k}(0x)\ =\ \alpha_{j,k}(x),$
\end{xlist}
so $\alpha_{j,k}(x)\in X^M_k,$ so the arbitrary element
$\alpha_{j,\infty}(x)\in (\coLm_I\,X_i)^M$
is in the image of~(\ref{iotaM}).

Recalling that an element $z$ of a monoid $M$ is called
a {\em right zero} element if $uz=z$ for all $u\in M,$
we see that the above argument shows that a sufficient condition
for~(\ref{bijM}) to hold, clearly independent of~(\ref{Mfg}), is
\begin{xlist}
\item\label{rtz}
$M$ has at least one right zero element.
\end{xlist}

With a little thought, one can come up with a common generalization
of~(\ref{Mfg}) and~(\ref{rtz}).
Recall that a {\em left ideal} of a monoid means a subset $L$
closed under left multiplication by all elements of $M.$
Combining the ideas of the two preceding arguments, one can show
that~(\ref{bijM}) holds if
\begin{xlist}
\item\label{fgI}
$M$ has a nonempty left ideal $L$ which is finitely generated
as a semigroup.
\end{xlist}

But we can generalize this still further.
We don't need left multiplication by {\em every} element of $M$ to send
{\em every} element of $L$ into $L.$
We claim it suffices to assume that
\begin{xlist}\item[]
$M$ has a finitely generated subsemigroup $S$ such
that $\{a\in M~|\linebreak[0] aS\cap S\neq\varnothing\}$ generates $M.$
\end{xlist}
Indeed, assuming the above holds, and given as before a directed system
$(X_i)_{i\in I}$ of $\!M\!$-sets and an element
$\alpha_{j,\infty}(x)\in (\coLm_I\,X_i)^M,$ choose
$k\geq j$ such that for all elements $g$ of a finite
generating set for $S,$
we have $g\,\alpha_{j,k}(x)=\alpha_{j,k}(x);$
thus $\alpha_{j,k}(x)$ is invariant under the action of $S.$
Writing $\alpha_{j,k}(x)=y,$ note that
for any $a\in M$ such that $aS\cap S\neq\varnothing,$
if we take $s,t\in S$ such that $as=t,$ and apply the two sides of
this equation to $y,$ we get $ay=y,$ showing that $y$ is fixed
under the action of each such element $a.$
Since such elements generate $M,$ we can conclude
that $y\in X^M_k,$ from which~(\ref{bijM}) follows as before.

In the condition just considered, nothing is lost if we replace
the semigroup
$S$ by the monoid $S\cup\{1\}.$
(The same was not true of~(\ref{fgI}), where the property of being
an ideal would have been lost.)
So let us formulate that condition in the more natural form
\begin{xlist}\item\label{fgM0}
$M$ has a finitely generated submonoid $M_0$ such
that $\{a\in M~|\linebreak[2] aM_0\cap M_0\neq\varnothing\}$
generates $M.$
\end{xlist}

To see that this is strictly weaker than~(\ref{fgI}), consider
the monoid presented by infinitely many generators
$x_n$ $(n\in\Nset)$ and $y,$
and the relations saying that all the elements
$x_n\<y~(n\in \Nset)$ are equal.
Then~(\ref{fgM0}) holds with $M_0$ the submonoid generated
by $\{y,\,x_0\<y\},$ but one can verify that there is no left
ideal $L$ as in~(\ref{fgI}).
(In particular, the left ideal $M\!\<y$ is {\em not} finitely
generated as a semigroup: the infinitely many elements
$x_n\<x_0\<y~(n\in \Nset)$ cannot be obtained using finitely
many elements of that ideal.)

Note that in condition~(\ref{fgM0}), one obtains the elements of $M_0$
from a finite generating set using arbitrarily many multiplications;
then gets each element $a$ in the set-bracket expression
from two elements of $M_0$ by an operation of ``right division'',
and then obtains the general element of $M$
from these by again using arbitrarily many multiplications.
Looked at this way, it would be more natural to allow arbitrary
sequences of multiplications and right divisions; i.e., to consider
the condition

\begin{xlist}\item\label{multdiv}
There exists a finite subset $S\subseteq M$ such that
the least subset $N\subseteq M$ satisfying
~(i)~$S\,\cup\,\{1\} \subseteq N,$
~(ii)~$a,b\in N\implies ab\in N$ and
~(iii)~$ab,b\in N\implies a\in N,$ is $M$ itself.
\end{xlist}

We shall see in the next section that this, too implies~(\ref{bijM}).
That~(\ref{multdiv}) is weaker than~(\ref{fgM0}) may be
seen by considering the monoid with presentation
\begin{xlist}\item[]
$M\ =\ \langle x_n,\,y_n,\,z,\,w\ \,(n\in \Nset)
\ |\ x_n\<y_n\<z=z,\ y_n\<w=w\rangle.$
\end{xlist}
Namely, one can show that
given a finitely generated submonoid $M_0\subseteq M,$
only finitely many of the elements $x_n$ can satisfy
$x_nM_0\cap M_0\neq\varnothing,$ hence not all $x_n$ will appear
in the set-expression shown in~(\ref{fgM0}), so, as these elements
are irreducible,~(\ref{fgM0}) cannot hold.
However, starting with the finite set $\{z,w\},$ the ``right
division'' process of~(\ref{multdiv}) gives us all elements of the
forms $x_n\<y_n$ and $y_n,$ another application
of right division gives all elements $x_n,$ and from the $y_n,$
the $x_n,$ and the original two elements $z$ and $w,$
closure under multiplication produces all of $M.$

\section{Left congruences, and a precise criterion.}\label{moncrit}

To approach more systematically the problem of characterizing
monoids that satisfy~(\ref{bijM}), let us recall a useful heuristic
for generalizing results about groups $G$ and $\!G\!$-sets to
monoids $M$ and $\!M\!$-sets:
\begin{xlist}\item\label{analogy}
Groups~~{\large:}~~normal subgroups~~{\large:}~~subgroups~~%
{\large:\::}~~\\%
monoids~~{\large:}~~congruences~~{\large:}~~left congruences.
\end{xlist}
Normal subgroups $N$ of a group $G$ classify the
homomorphic images $f(G)$ of $G,$ by listing the elements that
fall together with $1$ under $f.$
To determine the structure of a homomorphic image $f(M)$ of
a monoid $M,$ it is not sufficient to consider elements that fall
together with $1;$ instead one must look at the set of
{\em all} pairs of elements that fall together,
$C=\{(a,b)\in M{\times}M~|~ f(a)=f(b)\}.$
Sets $C$ that arise in this way are called {\em congruences} on $M;$
these are precisely the subsets $C\subseteq M{\times}M$ such that
\begin{xlist}\item\label{cong}
$C$ is an equivalence relation which is closed under left
and right translation by elements of~$M.$
\end{xlist}

When we study the structures of left $\!G\!$-sets $X$ for $G$ a group,
the key concept is the set $G_x$ of elements of $G$ fixing
a given $x\in X,$ which may be any subgroup.
For $M$ a monoid and $x$ an element of a left $M\!$-set,
the analogous entity is the set $C_x=\{(a,b)\in M{\times}M~|~ax=bx\}.$
This can be any subset $C\subseteq M{\times}M$ satisfying
\begin{xlist}
\item\label{lcong}
$C$ is an equivalence relation closed under left
translation by all elements of $M.$
\end{xlist}
Such a set is called a {\em left congruence} on $M.$

For $G$ a group, every $\!G\!$-set is a disjoint
union of orbits $Gx\cong G/H.$
There is no such simple structure theorem
for a set $X$ on which a monoid $M$ acts.
Nevertheless,
such an $X$ is, of course, a {\em union} of orbits $Mx\cong M/C_x,$
and this fact will allow us to reduce~(\ref{bijM})
to a condition on left congruences.

(Aside:  We have mentioned 2-sided congruences, i.e., sets
satisfying~(\ref{cong}), only for perspective.
{\em Right actions} of monoids lead to a third concept,
that of a right congruence, left-right dual to~(\ref{lcong}).
But since right actions of $M$ are
equivalent to left actions of the opposite monoid, we lose no generality
by restricting attention in this note to left $\!M\!$-sets.)

Given a monoid $M$ and a subset $R\subseteq M{\times}M,$
there is a least left congruence $C$ containing $R,$
the left congruence {\em generated} by $R,$ obtained
by closing $R$ under the obvious operations (one each to obtain
reflexivity, symmetry, transitivity, and left
translation by each element of $M).$
Thus, one can speak of a left congruence being {\em finitely generated}.

The whole set $M{\times}M$ constitutes the {\em improper}
left congruence on $M.$
We shall now show that the necessary and sufficient condition on a
monoid $M$ for~(\ref{bijM}) to hold is
\begin{xlist}
\item\label{Cfg}
The improper left congruence on $M$ is finitely generated.
\end{xlist}
Moreover, we will find that the final condition~(\ref{multdiv})
of the preceding section is also equivalent to this.

The reader who is inclined to skip the proof below as straightforward
should note that the step~(\ref{bijM})\!$\implies$\!(\ref{Cfg})
involves an unexpected hiccup; I therefore recommend reading at
least that step.

\begin{thm}\label{TMfixed}
If $M$ is a monoid, $I$ a directed partially ordered
set, and\linebreak[4] $(X_i,\alpha_{i,j})_{i,j\in I}$ a directed
system of $\!M\!$-sets, then the set-map $\iota$ of~{\rm(\ref{iotaM})}
is one-to-one.

Moreover, for any monoid $M,$ the following implications hold among
the conditions introduced above:
$$\begin{array}{rcr}
\vspace{.4em} & \mathrm{(\ref{Mfg})} & \\
\vspace{.3em} & \big\Downarrow & \\
\mathrm{(\ref{rtz})}\implies\!&\mathrm{(\ref{fgI})}&\!\implies
\mathrm{(\ref{fgM0})\implies
(\ref{multdiv})\iff(\ref{bijM})\iff(\ref{Cfg})}.
\end{array}$$
\end{thm}\begin{pf}
The first assertion and the
implications through~(\ref{multdiv}) have already been noted.
(Moreover, none of those implications is reversible;
examples were given where this was not obvious.)
We shall complete the proof by showing $\mathrm{(\ref{multdiv})\implies
(\ref{bijM})\implies (\ref{Cfg}) \implies (\ref{multdiv}).}$

Given a finite set $S$ as in~(\ref{multdiv}) and an element
$\alpha_{j,\infty}(x)\in(\coLm_I\,X_i)^M,$ let us
take $k\in I$ such that the finitely many relations
$s\,\alpha_{j,k}(x)=\alpha_{j,k}(x)~(s\in S)$ all hold, and
let $y=\alpha_{j,k}(x).$
Then it is easy to check that the set $N=\{s\in M~|\ sy=y\}$ satisfies
conditions~(i)-(iii) of~(\ref{multdiv}), hence is all of $M.$
Thus $y$ is an element of $X^M_k$ mapping to the
given element $\alpha_{j,\infty}(x)$ of $(\coLm_I\,X_i)^M,$
proving~(\ref{bijM}).

The proof that~(\ref{bijM})\!$\implies$\!(\ref{Cfg})
starts like the corresponding
argument for groups:  If the improper left congruence on $M$ is not
finitely generated, let $I$ be the set of all finitely
generated left congruences on $M,$ partially ordered by inclusion.
The $\!M\!$-sets $X_C=M/C~(C\in I)$ will form a directed system such
that $\coLm_I\,X_C$ is the $\!1\!$-element $\!M\!$-set; hence
$(\coLm_I\,X_C)^M\neq\varnothing;$ but we claim that
each set $X^M_C~(C\in I)$ is empty.

For assume, on the contrary, that $X^M_C$ were nonempty.
If $M$ were a group, that would make $X_C$ a singleton,
hence it would make $C$ the improper left congruence, a contradiction.
For $M$ a general monoid, we can only conclude that {\em some}
equivalence class $[a]\in X_C$ is fixed under the action of $M.$
However, given such an $[a],$ let $C'$ be the left congruence
on $M$ generated by $C$ and the one additional pair $(a,1).$
Then in $M/C'$ the generating element $[1]=[a]$ is $\!M\!$-fixed,
so $C'$ is the improper left congruence, this time indeed contradicting
the assumption that the latter is not finitely generated.

Finally, to show~(\ref{Cfg})\!$\implies$\!(\ref{multdiv}), suppose
$\{(a_1,b_1),\dots,(a_n,b_n)\}$ is a finite generating set for
the improper left congruence on $M.$
Let $S=\{a_1,\dots,a_n,\:b_1,\dots,\linebreak[0]b_n\},$ let $N$
be the set constructed from $S$ as in~(\ref{multdiv}),
and let $U\subseteq M{\times}M$ be the set of ordered pairs which
can be written $(as,at)$ with $a\in M$ and $s,t\in N.$
By the closure properties of $N$ we see that each
$(as,at)\in U$ either has both components in $N$ (if $a\in N,$
by~(\ref{multdiv})(ii)) or neither
(if $a\notin N,$ by~(\ref{multdiv})(iii)).
It follows that the least equivalence relation $C$ containing
$U$ will not relate elements in $N$ with elements not in $N.$
Moreover, $U$ is closed under left translation by members
of $M,$ hence so is $C,$ i.e., $C$ is a left congruence on $M.$
But $C$ contains $\{(a_1,b_1),\dots,(a_n,b_n)\},$ so by choice of this
set, $C$ must be the improper left congruence; hence as it does not
relate elements in $N$ with elements not in $N,$
we must have $N=M,$ establishing~(\ref{multdiv}).\qed
\end{pf}

We remark that none of conditions of the above theorem
except~(\ref{Mfg}) is right-left symmetric.
Indeed, let $M$ consist of the identity element and
an infinite set $S$ of right-zero elements.
Then $M$ satisfies~(\ref{rtz}), hence satisfies all these
conditions other than~(\ref{Mfg}), but I claim that the opposite
monoid $M^\mathrm{op}$ does not satisfy~(\ref{Cfg}), hence
does not satisfy any of the conditions shown.
For any equivalence relation on the underlying set
of a monoid respects both left multiplication
by the identity and left multiplication by any left zero element;
hence every equivalence relation on $M^\mathrm{op}$ is a
left congruence; but the improper equivalence relation on an infinite
set is not finitely generated.

Incidentally, there is a simpler example for monoids than for
groups showing that~(\ref{bijM}) can fail
but the analogous statement~(\ref{bijM}$\!_{\Nset}\!$)
on direct limits indexed by the
natural numbers hold; equivalently, that the improper left
congruence may be non-finitely generated, yet not expressible
as the union of a countable chain of proper left congruences.
Let $M=\omega_1,$ the first uncountable ordinal, made a monoid
under the commutative binary operation $\mathrm{sup}.$
Every left congruence on $M$ corresponds to a decomposition
into disjoint convex sets (i.e., intervals); let us associate to
each proper left congruence $C$ the least $\alpha\in\omega_1$
such that $(0,\alpha)\notin C.$
By considering the sequence of ordinals associated in this
way with a countable ascending chain of such left congruences,
we see that its union cannot be the improper left congruence.

Before leaving the topic of monoids and their left congruences, let me
mention a tantalizing open question of Hotzel~\cite{Hotzel}
(slightly restated): If a monoid $M$ has {\em ascending chain
condition} on left congruences, must $M$ be finitely generated?
An affirmative answer has been proved under the assumption of ascending
chain condition on {\em both} right and left
congruences~\cite{Kozhukhov}.
For some further observations see~\cite[Problem~1]{WTq}.

\section{Functors on small categories.}\label{cat}
As noted in the introduction, a monoid $M$ can be regarded as the
system of morphisms of a one-object category $\fb E.$
An $\!M\!$-set $X$ is
then equivalent to a functor $\fb E\rightarrow \fb{Set},$ and the
fixed-point set of the action of $M$
on $X$ is the {\em limit} of that functor.
In the remaining sections, we shall extend the ideas
of the preceding section by replacing fixed-point sets of monoid
actions with limits of set-valued
functors on a general small category.

If $\fb E$ is a small category we shall, to maintain
parallelism with preceding sections, call a covariant
functor $\fb E\to\fb{Set}$ an ``$\!\fb E\!\<$-set'', and
denote such functors by $X$ and neighboring letters.
Objects of $\fb E$ will generally be denoted $E,F,\dots$
and morphisms of $\fb E$ by letters $a,b,\dots\,.$
For $E, F\in\ObE,$ the set of morphisms $E\to F$ will be
written $\fb E(E,F).$ We will assume that
$\fb E(E,F)$ and $\fb E(E',F')$ are disjoint unless $E=E'$ and $F=F'.$
If $\alpha: X\to X'$ is a morphism of $\!\fb E\!\<$-sets, its component
set-maps will be denoted $\alpha(E): X(E)\to X'(E)$ $(E\in\ObE).$

We recall that if $X$ is an $\!\fb E\!\<$-set, then $\Lm_{\fb E}\,X$
can be constructed as the set of $\!\ObE\!$-tuples $x=(x_E)_{E\in\ObE},$
with $x_E\in X(E)$ for each $E\in\ObE,$ which satisfy
the ``compatibility'' conditions
\begin{xlist}
\item\label{compat}
$(\<\forall\,E,F\in\ObE,\,a\in\fb E(E,F))\ \ X(a)(x_E)\ =\ x_F.$
\end{xlist}

By a {\em directed system} of $\!\fb E\!\<$-sets we shall mean a
family of $\!\fb E\!\<$-sets $(X_i)_{i\in I}$ indexed by a nonempty
directed partially ordered set $I,$ and given with morphisms
of $\!\fb E\!\<$-sets $\alpha_{i,j}\<{:}\ X_i\to X_j~(i\leq j\in I)$
such that each $\alpha_{i,i}$ is the identity morphism
of the $\!\fb E\!\<$-set $X_i,$ and for $i\leq j\leq k\in I,$
one has $\alpha_{i,k}=\alpha_{j,k}\,\alpha_{i,j}.$

Given such a system, we see that for each $E\in\ObE,$ the sets
$X_i(E)$ $(i\in I)$ and set-maps
$\alpha_{i,j}(E)\<{:}\ X_i(E)\to X_j(E)$ form a directed system of sets.
If we take the direct limit of each of these systems,
functoriality of the direct limit construction yields,
for each morphism $a\in\fb E(E,F),$
a set-map $\coLm_{i\in I}\,X_i(E)\to\coLm_{i\in I}\,X_i(F)$
which we shall write $(\coLm_{i\in I}\,X_i)(a),$ and whose action
on elements is described by
\begin{xlist}\item\label{Lma}
$(\coLm_{i\in I}\,X_i)(a)\:
(\alpha_{j,\infty}(E)(y))\ =\ \alpha_{j,\infty}(F)(X_j(a)(y))\quad
(y\in X_j(E)).$
\end{xlist}
These maps together make the family of direct-limit sets
$(\coLm_{i\in I}\,X_i(E))_{E\in\ObE}$ into
an $\!\fb E\!\<$-set, which we shall denote $\coLm_{i\in I}\,X_i.$
(It is not hard to show that this $\!\fb E\!\<$-set is in fact the
direct limit, i.e., colimit \mbox{\cite[p.67]{CW}},
\mbox{\cite[\S\S7.5-7.6]{245}},
of the directed system $(X_i)_{i\in I}$ in the category of
$\!\fb E\!\<$-sets, though we shall not need that fact.)
As with any $\!\fb E\!\<$-set, we can take its
category-theoretic limit, getting a set \begin{xlist}\item[]
$\Lm_{\fb E}\,(\coLm_{i\in I}\,X_i).$
\end{xlist}

On the other hand, starting with our original directed system
$(X_i)_{i\in I}$ of $\!\fb E\!\<$-sets, we can take the limit
over $\fb E$ of each $\!\fb E\!\<$-set $X_i,$ getting a
system of sets $(\Lm_{\fb E}\,X_i)_{i\in I}.$
The functoriality of this limit construction
yields connecting maps which we may denote
\begin{xlist}\item[]
$\Lm_{\fb E}\,\alpha_{i,j}\,{:}\ \ \Lm_{\fb E}\,X_i\to
\Lm_{\fb E}\,X_j\quad (i\leq j\in I),$
\end{xlist}
so we may form the direct limit of these sets, getting a set
\begin{xlist}\item[]
$\coLm_{i\in I}\,(\Lm_{\fb E}\,X_i).$
\end{xlist}

And once again there is a natural set-map connecting these
constructions,
\begin{xlist}
\item\label{iotaE}
$\iota:\ \coLm_{i\in I}\,(\Lm_{\fb E}\,X_i)\longrightarrow\,
\Lm_{\fb E}\,(\coLm_{i\in I}\,X_i).$
\end{xlist}

To describe $\iota$ explicitly, consider an element
of $\coLm_{i\in I}\<(\Lm_{\fb E}\,X_i),$ written as
$\alpha_{i,\infty}(x)$ for some $i\in I$ and $x\in\Lm_{\fb E}\,X_i.$
Since $x$ is an $\!\ObE\!$-tuple $(x_E)$
satisfying~(\ref{compat}), we can apply
$\alpha_{i,\infty}(E)$ to each component $x_E,$
getting an $\!\ObE\!$-tuple of elements of the sets
$\coLm_{i\in I}\,X_i(E)~(E\in\ObE).$
The compatibility conditions~(\ref{compat}) on the
components $x_E$ of the given element $(x_E)$ imply the
compatibility of the components of the resulting
family $(\alpha_{i,\infty}(E)(x_E))_{E\in\ObE},$ so that this
becomes an element of $\Lm_{\fb E}\,(\coLm_{i\in I}\,X_i),$
which is easily shown to be independent of the choice of
expression $\alpha_{i,\infty}(x)$ for our
given element of $\coLm_{i\in I}\<(\Lm_{\fb E}\,X_i).$

This time, however, even injectivity of $\iota$ is not automatic.
To obtain a criterion for it to hold, we will use a lemma on
partially ordered sets.
Recall that a subset $D$ of a partially ordered set $J$ is called
a {\em downset} (or ``order ideal'') if $s<t\in D\implies s\in D.$
We shall regard the set of downsets of any partially ordered set
as ordered by inclusion.
A partially ordered set is called {\em downward directed}
(the dual of ``directed'') if for any
two elements $u,\,v$ of the set, there is an element $w$
majorized by both of them.

\begin{lem}\label{Lmin}
Let $J$ be a partially ordered set.
Then the following conditions are equivalent:
\begin{xlist}
\item\label{abovefin}
There exists a finite subset $A\subseteq J$ such that every
element of $J$ majorizes at least one element of $A\,.$
\end{xlist}
\begin{xlist}
\item\label{goodmin}
$J$ has only finitely many minimal elements, and every
element of $J$ majorizes a minimal element.
\end{xlist}
\begin{xlist}
\item\label{capne}
Every set $S$ of nonempty downsets of $J$ which is nonempty
and downward directed under inclusion has nonempty intersection.
\end{xlist}
\end{lem}\begin{pf}
Clearly~(\ref{goodmin})\!$\implies$\!(\ref{abovefin}).
To show (\ref{abovefin})\!$\implies$\!(\ref{capne}),
let $A$ be as in~(\ref{abovefin}), let $S$ be as in the
hypothesis of~(\ref{capne}), and for each
$a\in A$ which does not belong to all the elements of $S,$
choose an element $s(a)\in S$ not containing $a.$
Since $A$ is finite and $S$ is downward directed, we can find some
$s\in S$ which is majorized by (i.e., is a subset of) all these
sets $s(a).$
Being a nonempty downset, $s$
must contain some element of $A$ by~(\ref{abovefin}),
and by choice of $s$ that element belongs to all
members of $S,$ proving (\ref{capne}).

Finally, assuming~(\ref{capne}) we will prove~(\ref{goodmin}).
On the one hand,~(\ref{capne}), applied to chains $S$ and combined
with Zorn's Lemma (used upside down)
shows that every nonempty downset contains a minimal nonempty downset,
which must be a singleton consisting of a minimal element;
hence every element of $J$ majorizes a minimal element.
Moreover, if the set of minimal elements were infinite, then the
set $S$ of cofinite subsets of that set would be a counterexample
to~(\ref{capne}); so there are indeed only finitely many minimal
elements.\qed
\end{pf}

We can now get a criterion for the injectivity of
the set-maps $\iota,$ and a little more.

\begin{prop}\label{PinjE}
If $\fb E$ is a small category,
the following conditions are equivalent:
\begin{xlist}
\item\label{injE}
For every directed partially ordered set $I$ and
directed system $(X_i,\alpha_{i,j})_{i,j\in I}$ of $\!\fb E\!\<$-sets,
the set-map $\iota$ of~{\rm(\ref{iotaE})} is one-to-one.
\end{xlist}
\begin{xlist}
\item\label{frfin}
There exists a finite family $A$ of objects of $\fb E$
such that every object of $\fb E$ admits a
morphism from one of the objects of $A.$
\end{xlist}

Moreover, condition~{\rm(\ref{frfin})} is also {\em necessary}
for the map $\iota$ to be {\em surjective} for all directed systems.
\end{prop}\begin{pf}
First, assume~(\ref{frfin}), and let us
be given two elements $\alpha_{j,\infty}(x)$ and
$\alpha_{j',\infty}(y)$ in $\coLm_{i\in I}\,(\Lm_{\fb E}\,X_i)$
(where $x=(x_E)_{E\in\ObE}$ and $y=(y_E)_{E\in\ObE}),$
having the same image in $\Lm_{\fb E}\,(\coLm_{i\in I}\,X_i).$
Thus, the images of these two $\!\ObE\!$-tuples
agree in each component $\coLm_{i\in I}\,X_i(E)$ $(E\in\ObE).$
By the directedness of $I$ we can find $k$ majorizing
both $j$ and $j'$ and such that for each of
the finitely many objects $E\in A,$ $\alpha_{j,k}(E)(x_E)$
and $\alpha_{j',k}(E)(y_E)$ coincide.
Now by assumption, every $F\in\ObE$ admits a morphism from one
of the objects $E\in A,$ so the conditions~(\ref{compat})
on the $\!\ObE\!$-tuples $\alpha_{j,k}(x_E)$ and
$\alpha_{j',k}(y_E)$ show that the $\!F\!$-components of these
tuples coincide as well.
Hence $\alpha_{j,k}(x)=\alpha_{j',k}(y);$ hence
$\alpha_{j,\infty}(x)=\alpha_{j',\infty}(y),$ proving~(\ref{injE}).

To get the converse,
let us define a preordering on $\ObE$ by writing $E\leq F$
if there exists a morphism from $E$ to $F,$ and let $J$
be the partially ordered set obtained by dividing $\ObE$ by the
equivalence relation ``$\!E\leq F\leq E$''.
If~(\ref{frfin}) fails, this says that $J$ does not
satisfy~(\ref{abovefin}), hence by the preceding lemma we can
find a downward directed set $S$ of nonempty downsets
in $J$ having empty intersection.
We shall now construct a directed system of $\!\fb E\!\<$-sets
indexed by the (upward) directed partially ordered set $S^\mathrm{op}.$

Given $s\in S,$ let us say that an object $E\in\ObE$ ``belongs
to'' $s$ if the equivalence class of $E$ in $J$ is a member of $s.$
For each $s\in S,$ we define an $\!\fb E\!\<$-set $X_s$ by letting
$X_s(E)$ be the two-element set $\{-1,+1\}$ if $E$
belongs to $s,$ and the one-element set $\{0\}$ otherwise.
Given a morphism $a\in\fb E(E,F),$ we let $X_s(a)$ be
the identity on $\{-1,+1\}$ if $E$ and $F$ both belong
to $s;$ as $s$ is a downset, the remaining possibilities all
have $F$ not belonging to $s,$ in which case we
let $X_s(a)$ be the unique map $X_s(E)\to X_s(F)=\{0\}.$

If $s\supseteq t$ are members of $S,$ then we define the
map $\alpha_{s,t}\<{:}\ X_s\to X_t$ to act as the identity at
objects $E\in\ObE$ belonging either to both $s$ and $t$
or to neither, and as the unique map $\{-1,+1\}\to\{0\}$ on
elements belonging to $s$ but not to $t;$ these maps clearly
make $(X_s,\alpha_{s,t})_{s,t\in S}$
a directed system indexed by $S^\mathrm{op}.$
Now because $S$ has empty intersection, we see that at each
$E\in\ObE,$ the sets $X_s(E)$ become singletons for sufficiently large
$s\in S^\mathrm{op},$ so $\coLm_{S^\mathrm{op}}\,X_s$ is
an $\!\fb E\!\<$-set all of whose components are singletons; hence
the set $\Lm_{\fb E}\,(\coLm_{S^\mathrm{op}}\,X_s)$ is a singleton.

On the other hand, for each $s\in S$ we can construct
(at least) two distinct elements of $\Lm_{\fb E}\,X_s;$ an element
$x^+$ which takes value $+1$ at every $E$ belonging to
$s$ (and, necessarily, value $0$ at all other $E),$ and an
element $x^-$ which takes value $-1$ at all $E$ belonging
to $s.$
The maps $\Lm_{\fb E}\:\alpha_{s,t}\,{:}\ \ \Lm_{\fb E}\:X_s\to
\Lm_{\fb E}\:X_t$ $(s\supseteq t)$
take $x^+$ to $x^+$ and $x^-$ to $x^-;$ thus we get
distinct elements $x^+$ and $x^-$ in
$\coLm_{S^\mathrm{op}}\,(\Lm_{\fb E}\,X_s).$
Hence the map~(\ref{iotaE}) cannot be one-to-one.

The final assertion of the proposition is proved by a construction
exactly like the preceding, with $\varnothing$ used
in place of $\{-1,+1\}.$
In that case we get $\coLm_{S^\mathrm{op}}\,(\Lm_{\fb E}\,X_s)$
empty, and $\Lm_{\fb E}\,(\coLm_{S^\mathrm{op}}\,X_s)$ again
a singleton, so that~(\ref{iotaE}) is not surjective.\qed
\end{pf}

To formulate a criterion for~(\ref{iotaE}) to be bijective for
all directed systems of $\!\fb E\!\<$-sets, let us
define a {\em congruence} $C$ on an $\!\fb E\!\<$-set $X$ to be a family
$(C_E)_{E\in\ObE},$ where each $C_E$ is an equivalence relation
on $X(E),$ and which is functorial, in the sense that
\begin{xlist}
\item\label{eqrelE}
$(s,t)\in C_E,\ a\in\fb E(E,F) \implies (X(a)(s),\:X(a)(t))\in C_F.$
\end{xlist}
If, more generally, we define a ``binary relation'' $R$
on an $\!\fb E\!\<$-set $X$ to mean a family
$R=(R_E)_{E\in\ObE},$ where each $R_E$ is
a binary relation on $X(E),$ and no functoriality is assumed,
then for every such relation $R$ we can define the
congruence {\em generated by} $R$ to be the least congruence
$C$ such that for each $E\in\ObE,$ $R_E\subseteq C_E.$
It is not hard to verify a more explicit description for this
congruence $C\,{:}$ for each $E\in\ObE,$
$C_E$ is the equivalence relation on $X(E)$ generated by
the union, over all $F\in\ObE$ and $a\in\fb E(F,E),$ of the image
in $X(E){\,\times\,}X(E)$ of $R_F\subseteq X(F){\,\times\,}X(F)$
under $X(a){\,\times\,}X(a).$
We will call a congruence on $X$ {\em finitely generated} if it
is generated by a binary relation $R$ such
that $\sum_{E\in\ObE} \mathrm{card}(R_E)<\infty.$
(Since we cannot assume the sets $X(E)$ disjoint, it is not
sufficient to say that $\mathrm{card}(\,\bigcup R_E)$ is finite.)
The {\em improper} congruence on an $\!\fb E\!\<$-set $X$ will mean
the congruence whose value at
each $E$ is the improper equivalence relation on $X(E).$

For any object $E$ of $\fb E,$ the covariant
hom-functor $\fb E(E,-)\,{:}\,\ \fb E\to\fb{Set}$ may be regarded as
an $\!\fb E\!\<$-set, which we will denote $H_E.$
Since we have assumed that distinct pairs of objects
have disjoint hom-sets, these $\!\fb E\!\<$-sets will be disjoint,
and we can form the union of any set of them.
We can now state and prove

\begin{thm}\label{TEfixed}
Let $\fb E$ be a small category satisfying~{\rm(\ref{frfin})},
and $A$ a finite set of objects of $\fb E$ as in that condition, i.e.,
such that every object of $\fb E$ admits at least one morphism from an
object of $A.$
Let $H$ denote the $\!\fb E\!\<$-set $\:\bigcup_{E\in A} H_E.$
Then the following conditions are equivalent:
\begin{xlist}
\item\label{bijE}
For every directed partially ordered set $I$ and
directed system $(X_i,\alpha_{i,j})_{i,j\in I}$ of $\!\fb E\!\<$-sets,
the set-map $\iota\<{:}\:\ \coLm_{i\in I}\,(\Lm_{\fb E}\,X_i)\to
\linebreak[4]\Lm_{\fb E}\,(\coLm_{i\in I}\,X_i)$
of~{\rm(\ref{iotaE})} is bijective.
\end{xlist}
\begin{xlist}
\item\label{Hfg}
The improper congruence on $H$ is finitely generated.
\end{xlist}
\end{thm}\begin{pf}
Since~(\ref{frfin}) is equivalent to injectivity
of the maps~(\ref{iotaE}), what we must prove is that under that
assumption, surjectivity of all such maps is equivalent to~(\ref{Hfg}).

First assume~(\ref{Hfg}), and suppose we are given a directed
system $(X_i)_{i\in I}$ of $\!\fb E\!\<$-sets, and an
element $x=(x_E)\in\Lm_{\fb E} (\coLm_{i\in I}\,X_i).$
Each coordinate $x_E$ of $x$ can be written
$\alpha_{j_E,\infty}(y_E),$ where $j_E\in I$ is an index
depending on $E,$ and $y_E\in X_{j_E}(E).$
We shall only use finitely many of these elements,
namely those with $E\in A.$
By the directedness of $I$ we can find an index $j$ that
majorizes all the $j_E$ with $E\in A;$ we thus get a family of elements
of $X_j,$ namely $y'_E=\alpha_{j_E,j}(y_E)\in X_j(E)$ $(E\in A).$
These will generate a sub-$\!\fb E\!\<$-set $Y\subseteq X_j,$
whose $\!F\!$-component, for each $F\in\ObE,$ consists of
all elements $X_j(a)(y'_E)$ $(E\in A,\:a\in\fb E(E,F)).$

Let us map the $\!\fb E\!\<$-set $H=\bigcup_{E\in A} H_E$
onto $Y$ by sending each $a\in H_E(F)=\fb E(E,F)$ (where
$E\in A,\:F\in\ObE)$ to $X_j(a)(y'_E)\in X_j(F).$
(This can be thought of as an application of Yoneda's Lemma to
each of the sub-$\!\fb E\!\<$-sets $H_E$ $(E\in A)$ of $H.)$

By choice of the $y_E,$ the image in $\coLm_{i\in I}\,X_i$
of the sub-$\!\fb E\!\<$-set $Y\subseteq X_j$
has in each coordinate $F$ only a single element, namely $x_F.$
Thus by applying the morphism $\alpha_{j,k}$ for large enough $k,$ we
can make any given pair of elements in any coordinate fall together.
But the fact that $Y$ is an image of $H$ and that the
improper congruence on $H$ is finitely generated means that
some finite family of these collapses imply all of them.
Thus, we can find some $k\geq j$ such that the image
of $Y$ in $X_k$ has just one element in each coordinate.
The $\!\ObE\!$-tuple of elements of $X_k$ so determined will
be an element $z\in\Lm_{\fb E}\,X_k$ which maps to $x$
in $\Lm_{\fb E} (\coLm_{i\in I}\,X_i).$
Taking the image of this element $z$ in
$\coLm_{i\in I}\,(\Lm_{\fb E}\,X_i)$
we get an element of the latter set that maps to
$x$ under $\iota,$ proving~(\ref{bijE}).

The proof of the converse will also follow that of
the corresponding result for monoid actions, though this time
the ``hiccup'' will involve adjoining $\mathrm{card}(A)$
additional pairs, rather than just one,
to a certain finitely generated congruence.
Assuming~(\ref{bijE}), let $I$ be the
directed partially ordered set of all finitely generated
congruences on $H,$ and for each $C\in I,$
let $X_C$ be the $\!\fb E\!\<$-set $H/C.$
Then we see that $\coLm_{C\in I}\,X_C$ is an $\!\fb E\!\<$-set
with just one element in each component, hence
$\Lm_{\fb E} (\coLm_{C\in I}\,X_C)$ is a singleton.
Hence by~(\ref{bijE}) the same is true of
$\coLm_{C\in I}\,(\Lm_{\fb E}\,X_C),$ so
least one of the sets $\Lm_{\fb E}\,X_C$ $(C\in I)$ is nonempty.
Say $x=(x_E)\in\Lm_{\fb E}\,X_C$ for some $C\in I.$
For each $E\in\ObE$ the element $x_E$ will be the
$\!C\!$-congruence class $[a_E]$ of some element
$a_E\in \fb E(F_E,E)\subseteq H(E),$ where $F_E\in A.$
If for every $F\in A$ we adjoin to $C$ the additional pair
$(\mathrm{id}_F,a_F),$ we get a congruence $C'$ on $H$ which is still
finitely generated, and which I claim is the improper congruence.
Indeed, the compatibility conditions~(\ref{compat}), which by
assumption hold for the components $x_E=[a_E]$ of $x,$ now hold
also for all translates $[a]$ $(a\in\fb E(F,E))$ of the
images $[\mathrm{id}_F]$ of the generators $\mathrm{id}_F$ of $H.$
This establishes~(\ref{Hfg}).\qed
\end{pf}

The following terminology provides a useful way of looking at this
result.

\begin{defn}\label{Dfinpres}
Let $\fb E$ be a small category.
By the {\em trivial $\!\fb E\!\<$-set} we will mean the
functor $T$ that associates to every object of $\fb E$ a $\!1\!$-element
set \r(\!with the only possible behavior on morphisms\r).
If $\fb E$ satisfies {\rm(\ref{frfin})} and \r(for $H$ constructed
as in Theorem~\ref{TEfixed}
from a set $A$ as in~{\rm(\ref{frfin}))}, also~{\rm(\ref{Hfg})},
we will say that the trivial $\!\fb E\!\<$-set is
{\em finitely presented}.
\end{defn}

Note that for $T$ the trivial $\!\fb E\!\<$-set defined above, and
$X$ any $\!\fb E\!\<$-set, the set $\Lm_{\fb E}\,X$ can be identified
with the hom-set $\fb{Set^E}(T,X).$
From this point of view, Theorem~\ref{TEfixed} is
an instance of the general observation that for an algebraic structure
$S$ (in this case, $T),$ the covariant hom-functor determined by
$S$ respects direct limits if and only if $S$ is finitely presented.

(We have, for simplicity, not defined the general concept of
a presentation of an $\!\fb E\!\<$-set.
Briefly, this may be done as follows.
A representable functor $H_E$ $(E\in\ObE)$
can be considered an $\!\fb E\!\<$-set $X$ ``free
on one generator in $X(E)\!$'', namely $\mathrm{id}_E.$
A disjoint union of $\!\fb E\!\<$-sets of this form (with repetitions
allowed), modulo the congruence generated by a given set of ordered
pairs, can be regarded as the $\!\fb E\!\<$-set presented using
the images of the elements $\mathrm{id}_E\in H_E$ as
generators and the given ordered pairs as relations.
Incidentally, the reader may have noted
that Definition~\ref{Dfinpres} has the formal defect
that the condition on $\fb E$ as stated depends on the choice of $A.$
But Theorem~\ref{TEfixed} shows that it is in fact independent
of $A;$ and, indeed, for $\!\fb E\!\<$-sets as for other finitary
algebraic objects, if an object is finitely generated, one can show
that the property of finite relatedness is independent of one's
choice of finite generating set.)

Though Theorem~\ref{TEfixed} is elegant, it does not give
convenient conditions analogous
\mbox{to~(\ref{Mfg})-(\ref{multdiv})} of Theorem~\ref{TMfixed}.
These, too, may be generalized to arbitrary small categories
$\fb E,$ but the statements are simplest when $\fb E$
has only finitely many objects.
I will develop the generalization of~(\ref{multdiv}) to that case below,
and at the end of the next section will state and sketch the proof of
the corresponding result for general $\fb E.$

Let us call a subcategory $\fb E_0$ of a category $\fb E$
{\em right division-closed} if for any two morphisms $a,b$ of $\fb E$
whose composite $ab$ is defined, we have
\begin{xlist}\item\label{divcl}
$ab,\<b\in\fb E_0\implies a\in\fb E_0.$
\end{xlist}

\begin{prop}\label{Pfinobjs}
Let $\fb E$ be a category with only finitely many objects.
Then the following conditions are equivalent:
\begin{xlist}
\item\label{cdnTEf}
$\fb E$ satisfies the equivalent conditions of Theorem~\ref{TEfixed}.
\end{xlist}
\begin{xlist}
\item\label{Emultdiv}
There exists a finite set $S$ of morphisms of $\fb E,$
such that the smallest subcategory $\fb E_0$
of $\fb E$ which has the same object-set as
$\fb E,$ and contains $S,$ and is
right division-closed in $\fb E,$ is $\fb E$ itself.
\end{xlist}
\end{prop}\begin{pf}
Assuming~(\ref{cdnTEf}), take for $A$ as in Theorem~\ref{TEfixed} the
full object-set of $\fb E,$ so that $H$ is the union of
the $\!\fb E\!\<$-sets $H_E$ associated with all the
objects of $\fb E,$ and let $R$ be a finite generating set
for the improper congruence on $H\!.$
Let $S$ be the set of all elements occurring as first
or second components of members of $R,$ and let $\fb E_0\subseteq \fb E$
be constructed from $S$ as in~(\ref{Emultdiv}).
Let $U$ be the set of all pairs $(as,at)$ with
$s\in\fb E_0(E,F),\ t\in\fb E_0(E',F),\ a\in\fb E(F,G),
\ E,E',F,G\in\ObE.$
As in the last paragraph of the
proof of Theorem~\ref{TMfixed}, each element of $U$
either has both components or neither component in $\fb E_0.$
Hence the equivalence relation $C$ generated by $U$ also has this
property.
Moreover, $U,$ and hence $C,$ is closed under left composition
with morphisms of $\fb E,$ hence $C$ is a congruence
on $H,$ and it contains all members of the generating set $R,$
hence it is the improper congruence.
Now for every morphism $a\in\fb E(E,F)$ of $\fb E,$ the improper
congruence on $H$ contains $(\mathrm{id}_F,a),$
and $\mathrm{id}_F\in\fb E_0,$ hence by our ``both or
neither'' property of $C,$ $a\in\fb E_0.$
So $\fb E_0=\fb E,$ proving~(\ref{Emultdiv}).

Conversely, suppose $S$ is a finite set of morphisms for which the
conclusion of~(\ref{Emultdiv}) holds, and consider the congruence $C$
on $H$ generated by all pairs $(a,\mathrm{id}_F)$ where
$a\in S\cap\fb E(E,F)$ $(E,F\in\ObE).$

For each $E,F\in\ObE,$ let $\fb E_1(E,F)$ denote
$\{a\in\fb E(E,F)~|~(a,\mathrm{id}_F)\in C\}.$
I claim this gives the morphism-set of
a right division closed subcategory
$\fb E_1\subseteq\fb E$ with object-set $\ObE.$
It is immediate that it contains all identity morphisms;
now suppose $a\<{:}\ F\to G$ and $b\<{:}\ E\to F$ are morphisms
of $\fb E,$ with $b\in\fb E_1.$
The latter relation means $(b,\mathrm{id}_F)\in C,$ hence as $C$ is
a congruence, we also have $(ab,a)\in C,$ hence
$(ab,\mathrm{id}_G)\in C{\iff}(a,\mathrm{id}_G)\in C,$
i.e., $ab\in\fb E_1{\iff}a\in\fb E_1,$ proving both closure
under composition and right division closure.
Hence since $S$ was chosen as
in~(\ref{Emultdiv}), $\fb E_1$ must be all of $\fb E.$
This says that $C$ contains all pairs
$(a,\mathrm{id}_F)$ with $a\in\fb E(E,F),\ E,F\in\ObE,$
so by transitivity, $C$ is the improper congruence
on $H,$ which is thus finitely generated, proving~(\ref{cdnTEf}).\qed
\end{pf}

\section{Digression: four corollaries.}\label{4cors}
A case of Proposition~\ref{Pfinobjs} which has been noted before is

\begin{cor}[{\rm $=$ \cite[Prop.~7.9.3]{245},
cf.\ \cite[Thm.~IX.2.1, p.211]{CW},
\cite[Thm.~4.73, p.72]{MK}}\textbf{}]
\label{CfgE}
If $\fb E$ is a category with only finitely many objects,
and whose morphism-set is finitely generated under composition,
then on directed systems $(X_i)_{i\in I}$ of $\!\fb E\!\<$-sets,
the operations $\Lm_{\fb E}$ and $\coLm_I$ commute;
i.e.,~{\rm(\ref{bijE})} holds.\qed
\end{cor}

This note is in fact the result of pondering how to improve on
the above result from \cite{245}.
(Incidentally, in the statement in \cite{245}, I assumed
the category $\fb E$ nonempty, but allowed direct limits over possibly
empty directed partially ordered sets.
In this note I have made the opposite choices, requiring in the
definition that direct limits have nonempty index sets, but
not so restricting $\fb E.$
As observed in \mbox{\cite[Exercise~7.9:2]{245}}, the result holds
if either the category or the directed partially ordered set is
required to be nonempty, but fails when they are both empty.)

Let us record next a pair of results implicit in the proofs of
Theorem~\ref{TMfixed} and Proposition~\ref{Pfinobjs}, along with
their duals.
(The proofs, and that of the final result of this section, will
just be sketched; they will not be used in the remainder of this note.)
We will use for monoids as well as for categories
the term ``right division-closed'' introduced above, and define
``left division-closed'' for both sorts of structures dually.
(The terms used by semigroup theorists, e.g. in \cite{C+P.2,BMS},
are ``left, respectively, right unital'', though
in~\cite{PD}, where the conditions were first introduced,
they were ``left, respectively, right simplifiable''.)
For $a$ an element of a monoid $M$ and $Y_0$ a subset of an
$\!M\!$-set $Y,$ we shall in Corollary~\ref{N=M_x}(iii*)
write $a^{-1}Y_0$ for the inverse image
of $Y_0\subseteq Y$ under the map $Y\rightarrow Y$ given by the
action of $a;$ and we shall similarly use inverse image notation
in Corollary~\ref{E_0=E_x}(iii*) in connection with the set-maps
$Y(a)$ forming the structure of an $\!\fb E\!\<$-set $Y.$
Note that in Corollary~\ref{N=M_x}(iii*), ``$\!M\!$-set''
still means left $\!M\!$-set, despite the dualization being carried out.

The first half of the next result is due to Schein
\cite[Theorem~2]{BMS}. 

\begin{cor}[{\rm to proof of Theorem~\ref{TMfixed};
cf.\ \cite{BMS}}\textbf{}]
\label{N=M_x}
Let $M$ be a monoid, and $N$ a subset of $M.$
Then the following conditions are equivalent:

{\rm(i)}~ $N$ is a right division-closed submonoid of $M.$

{\rm(ii)}~ $N$ is the equivalence class of $1$ under
some left congruence on $M.$

{\rm(iii)}~ There exist an $\!M\!$-set $X,$ and an element $x\in X,$
such that $N=\\ \{a\in M\ |\ ax=x\}.$

Likewise, the following conditions are equivalent:

{\rm(i*)}~ $N$ is a {\em left} division-closed submonoid of $M.$

{\rm(ii*)}~ $N$ is the equivalence class of $1$ under
some {\em right} congruence on $M.$

{\rm(iii*)}~ There exist an $\!M\!$-set $Y,$ and a subset
$Y_0$ of $Y,$ such that $N=\\ \{a\in M\ |\ a^{-1}Y_0=Y_0\}.$
\r(I.e., $N$ is the set of elements $a\in M$ which carry
both $Y_0$ and its complement into themselves.\r)
\end{cor}
\begin{pf*}{Sketch of proof.}
Assuming~(i), let $C$ be the equivalence relation
on $M$ generated by $\{(as,at)\ |\ a\in M,~ s,t\in N\},$ and
observe as in the proof of Theorem~\ref{TMfixed},
(\ref{Cfg})\!$\implies$\!(\ref{multdiv}), that $C$ is a congruence
and relates elements of $N$ only with elements in $N,$
establishing~(ii).
The implications (ii)\!$\implies$\!(iii)\!$\implies$\!(i) are
straightforward.

The second half of the result will follow from the first by left-right
dualization if we can establish that (iii*) is equivalent to the
existence of a {\em right} $\!M\!$-set $X$ with an element $x$
such that $N=\{a\in M\ |\ xa=x\}.$
Now given a left $\!M\!$-set $Y$ and a subset $Y_0$
such that $N=\{a\in M\ |\ a^{-1}Y_0=Y_0\},$ the
contravariant power functor yields a right $\!M\!$-set $X=\fb P(Y),$
in which the element $x=Y_0$ indeed satisfies
$N=\{a\in M\ |\ xa=x\}.$
Conversely, given a right $\!M\!$-set $X$
with an element $x$ satisfying this relation, it is easy to
verify that in the left $\!M\!$-set $Y=\fb P(X),$
the subset $Y_0 = \{S\subseteq X\ |\ x\in S\}$
satisfies $N=\{a\in M\ |\ a^{-1}Y_0=Y_0\}.$
\qed
\end{pf*}

We can now see the significance of condition~(\ref{multdiv})
in Theorem~\ref{TMfixed}.
Although, as noted at the beginning of \S\ref{moncrit},
a general left congruence on a monoid is not
determined by the set of elements congruent to~1, the improper left
congruence is clearly determined by that set.
Condition~(\ref{multdiv}) translates finite generation of the
improper left congruence into finite generation of $M$ as a set that
can occur as the equivalence class of $1$ under a left congruence.

The analogous result for small categories is

\begin{cor}[{\rm to proof of Proposition~\ref{Pfinobjs}}\textbf{}]
\label{E_0=E_x}
Let $\fb E$ be a small category, and for every pair of objects
$E,F\in\ObE$ let $\fb E_0(E,F)$ be a subset of $\fb E(E,F).$
Then the following conditions are equivalent:

{\rm(i)}~ The sets $\fb E_0(E,F)$ are the morphism-sets of a
right division-closed subcategory $\fb E_0\subseteq\fb E$ with the
same object-set as $\fb E.$

{\rm(ii)}~ There exists a congruence $C$ on the $\!\fb E\!$-set
$\bigcup_{E\in\ObE}H_E,$ such that for all
$E,F\in\ObE,$ $\fb E_0(E,F)=\{a\in \fb E(E,F)~|\ (a,1_F)\in C\}.$

{\rm(iii)}~ There exist an $\!E\!\<$-set $X,$ and for each
$E\in\ObE$ an element $x_E\in X(E),$ such that for all
$E,F\in\ObE,$ $\fb E_0(E,F) =\{a\in\fb E(E,F)\ |\ ax_E=x_F\}.$

Likewise, the following conditions are equivalent \r(\!\<where for
$E\in\ObE,$ $H^E$ denotes the contravariant hom functor
$\fb E(-,E)):$

{\rm(i*)}~ The sets $\fb E_0(E,F)$ are the morphism-sets of a
{\em left} division-closed subcategory $\fb E_0\subseteq\fb E$ with the
same object-set as $\fb E.$

{\rm(ii*)}~ There exists a congruence $C$ on the {\em right}
$\!\fb E\!$-set \r(contravariant $\!\fb{Set}\!$-valued functor\r)
$\bigcup_{E\in\ObE}H^E$ such that for all
$E,F\in\ObE,$ $\fb E_0(E,F)=\{a\in \fb E(E,F)~|\ (a,1_E)\in C\}.$

{\rm(iii*)}~ There exist a \r(left\/\r) $\fb E$-set $Y,$
and for each $E\in\ObE$ a subset $Y_0(E)\subseteq Y(E),$
such that for all $E,F\in\ObE,$
$\fb E_0(E,F) =\{a\in\fb E(E,F)\ |\\ \ Y(a)^{-1}Y_0(F) = Y_0(E)\}.$
\end{cor}
\begin{pf*}{Sketch of proof.}
Analogous to the proof of Corollary~\ref{N=M_x}.
So, for instance, to get (i)\!$\implies$\!(ii), we use $\fb E_0$
as in the proof of Proposition~\ref{Pfinobjs} to construct on $H$
a binary relation $U,$ and from that, the congruence $C.$\qed
\end{pf*}

Since this relation between right division-closed subcategories
and congruences on $H$ holds for arbitrary $\fb E,$
why does Theorem~\ref{TEfixed} need the hypothesis
that $\fb E$ have only finitely many objects?
Because when it has infinitely many
objects, the $\!\fb E\!$-set $\bigcup_{E\in\ObE}H_E$ is not
finitely generated, so the statement that its quotient by the
improper congruence is finitely presented does not mean that the
latter congruence is finitely generated.
However, with this viewpoint in mind, one can come up with
a generalization of that theorem to arbitrary $\fb E.$

\begin{cor}[{\rm to Theorem~\ref{TEfixed} and
proof of Proposition~\ref{Pfinobjs}}\textbf{}]
\label{S0S1}
Let $\fb E$ be a small category satisfying~{\rm(\ref{frfin})},
and $A$ a finite set of objects of $\fb E$ as in that condition.
Let $S_0$ be a set of morphisms of $\fb E$ which, for each
$F\in\ObE-A,$ contains exactly one morphism from a member
of $A$ to $F,$ and which contains no elements other than these.
Then $\fb E$ satisfies the equivalent conditions of
Theorem~\ref{TEfixed} if and only if it satisfies
\begin{xlist}
\item\label{01multdiv}
There exists a finite set $S_1$ of morphisms of $\fb E$
such that the smallest subcategory $\fb E_0$
of $\fb E$ which has the same object-set as
$\fb E,$ and contains $S_0\cup S_1,$ and is
right division-closed in $\fb E,$ is $\fb E$ itself.
\end{xlist}
\end{cor}
\begin{pf*}{Sketch of proof.}
(\ref{01multdiv}) is equivalent to the statement that
the pairs $(a,\mathrm{id}_F),$ where $a\in S_0\cup S_1$ and $F$
is the codomain of $a,$ generate the improper congruence
on $\bigcup_{E\in\ObE}H_E.$
Now those pairs with $a$ taken from $S_0$ simply
serve to ``eliminate'' the generators $\mathrm{id}_F$
of $\bigcup_{E\in\ObE}H_E$ with $F\in \ObE-A;$ i.e.,
dividing out by the congruence generated by those pairs alone
gives the \!$\fb E$\!-set $H$ of Theorem~\ref{TEfixed}.
Thus,~(\ref{01multdiv}) is equivalent to the statement that
the improper congruence on {\em that} \!$\fb E$\!-set is
generated by a finite set of pairs, which
is the desired condition~(\ref{Hfg}).

To set up a formal proof, for each $F\in\ObE-A,$ let
$a_F\in\fb E(E_F,F)$ (where $E_F\in A)$ be the corresponding
element of $S_0,$ while for $F\in A$ let us set $E_F=F,\ a_F=
\mathrm{id}_F.$
Then given $S_1$ as in~(\ref{01multdiv}), one shows that
the improper congruence on $H$ is generated by the finite set of pairs
$(b a_F,a_{F'})$ for $b:F\rightarrow F'$ in $S_1,$ while conversely,
given a finite generating set $R$ for that improper congruence, one
can take $S_1$ to be the set of components of members of $R\<.$\qed
\end{pf*}

An easy class of examples are categories $\fb E$ having
an initial object $E_\mathrm{init}.$
Then if one takes $A=\{E_\mathrm{init}\},$ there is a
unique set $S_0$ as in the statement of the above theorem, and
letting $S_1$ be the empty set, one finds that~(\ref{01multdiv}) holds.

\section{Posets.}\label{poset}

Groups and monoids, with which we began this note, are
categories where ``all the structure is in the morphisms'', and
essentially
none in the class of objects and the way morphisms connect them.
In this section we will consider the opposite extreme, the case of
partially ordered sets $J$ regarded as categories.

If $J$ is a poset,
we shall write $\fb E = J_\fb{cat}$ for the category having
for objects the elements of $J,$ and having, for
each $E, F\in J,$ one morphism $\lambda(E,F)\,{:}\ E\to F$ if $E\leq F,$
and no morphisms $E\to F$ otherwise.
(We write $E,~F,\dots$ for elements of $J$ for
consistency with the notation of the last two sections.)

From Proposition~\ref{PinjE} we know that a necessary condition for
limits over such a category $\fb E$ to respect direct limits is
that the set $A$ of minimal elements of $J$ be finite,
and every element of $J$ lie above an element of $A.$
Note that the $\!\fb E\!\<$-set $H$ constructed as in
Theorem~\ref{TEfixed}
from this set $A$ associates to each $E\in J$ the set
$\{\lambda(F,E)~|~F\,{\in}\<A,\,F\,{\leq}\,E\}.$
By that theorem, to strengthen our necessary condition to a
necessary and sufficient one,
we need to know for which $J$ the improper congruence on this
$\!\fb E\!\<$-set $H$ is finitely generated.

For an instructive example of an infinite poset for which
this congruence {\em is} finitely generated, let the underlying
set of $J$ consist of all real numbers $\geq 1,$ ordered in
the usual way, together with two elements $0_1$ and $0_2$ which
are less than all other elements, and mutually incomparable.
Thus, $A=\{0_1,0_2\},$ and for all $E$ other than these
two elements, we have $H(E)=\{\lambda(0_1,E),\lambda(0_2,E)\}.$
It is easy to see that the improper congruence
on $H$ is generated by the single pair
$(\lambda(0_1,1), \lambda(0_2,1)).$
On the other hand, if we delete the element $1$ and consider the
corresponding functor on $(J-\{1\})_\fb{cat},$ it is not hard to see
that the improper congruence on this
functor is no longer finitely generated.

The element $1\in J$ is what we shall call a ``critical element''
with respect to the subset $\{0_1,0_2\}.$
In the example above, it served to ``gather'' the strands
of $H$ emanating from $0_1$ and $0_2.$
Let us give precise meanings to these terms.

\begin{defn}\label{Dgath}
Let $J$ be a partially ordered set.
For $E\in J$ we shall write $\D(E)$ for $\{F\in J~|~
F\leq E\}$ \r(\!\<the ``principal downset'' generated by $E).$

Given $E\in J$ and subsets $A,B\subseteq J,$
we shall write $R(A,B,E)$ for the equivalence relation
on $A\<\cap~\D(E)$ generated by the union
over all $F\in B\cap \D(E)$ of the improper equivalence
relations on the sets $A\cap \D(F).$
We shall say that $B$ {\em gathers $A$ under $E$}
if $R(A,B,E)$ is the improper equivalence relation on
$A\cap \D(E).$

Given a subset $A\subseteq J$ and an element $E\in J,$
we note that $\{E\}$ always gathers $A$ under $E.$
We shall call $E$ {\em $A\!$-critical} if
$J-\{E\}$ does {\em not} gather $A$ under $E.$
\end{defn}

It is straightforward to verify the transitivity relation
\begin{xlist}
\item\label{trans}
If $A,\ B_1,\ B_2$ are subsets of $J$ and $E$ an
element of $J,$ such that $B_1$ gathers $A$ under
every element of $B_2$ and $B_2$ gathers $A$ under
$E,$ then $B_1$ gathers $A$ under $E.$
\end{xlist}
Also, the next-to-last sentence of
Definition~\ref{Dgath} implies the reflexivity condition:
\begin{xlist}
\item\label{refl}
If $A,\:B$ are subsets of $J,$ then $B$
gathers $A$ under every $E\in B.$
\end{xlist}

Note that in the next lemma, we do not assume that
every element of $J$ majorizes some member of $A$ (though
we will add that assumption when we apply the lemma).

\begin{lem}\label{Lgath}
Let $J$ be a partially ordered set and $A\subseteq J$ a finite
subset.
Let us write $\fb E$ for $J_\fb{cat},$ and $H$ for the
union, over all $E\in A,$ of the covariant hom-functors $H_E.$
Then the following conditions are equivalent:
\begin{xlist}
\item\label{finB}
There exists a finite subset $B\subseteq J$ which
gathers $A$ under every $E\in J.$
\end{xlist}
\begin{xlist}
\item\label{goodcrit}
The set of $A\!$-critical elements of $J$ is finite,
and gathers $A$ under every $E\in J.$
\end{xlist}
\begin{xlist}
\item\label{HfgJ}
The improper congruence on the $\!\fb E\!\<$-set $H$
is finitely generated.
\end{xlist}
\end{lem}\begin{pf}
(\ref{goodcrit})\!$\implies$\!(\ref{finB}) is immediate.
To get the converse, take $B$ as in~(\ref{finB})
and let $B'$ denote the set of $A\!$-critical elements of $J.$
Applying~(\ref{finB}) to an element $E\in B',$ we see, from the
definition
of the statement that $E$ is $\!A\!$-critical, that $E\in B.$
Hence $B'\subseteq B,$ so in particular $B'$ is finite; it remains
to show that for any $E\in J,$ $B'$ gathers $A$ under $E.$
In doing so we may assume inductively that $B'$ gathers $A$
under every $F\in J$ such that the number of elements
of $B$ that are $<F$ is smaller than the number that are $<E,$
or such that these numbers are equal
but the number $\leq F$ is smaller than the number $\leq E.$

If $E\notin B,$ the former assumption shows
that $B'$ gathers $A$ under each element of $B\cap \D(E),$
hence~(\ref{trans}), with $B'$ and
$B\cap \D(E)$ in the roles of $B_1$ and $B_2,$ shows
that $B'$ gathers $A$ under $E,$ as desired.
On the other hand, if $E\in B,$ the inductive assumptions
show that $B'$ gathers $A$ under every element of $J$ that is $<E.$
Now if $E$ is not $\!A\!$-critical, we can apply~(\ref{trans})
with $B'$ and $\D(E)-\{E\}$ in the roles of $B_1$ and $B_2$
respectively, and again conclude that $B'$ gathers $A$ under $E.$
On the other hand, if $E$ is $\!A\!$-critical, then it belongs
to $B',$ and~(\ref{refl}) (with $B'$ in the
role of $B)$ yields the same conclusion.

(\ref{finB})\!$\iff$\!(\ref{HfgJ}):
Note that for any $E\in J,$ the definition of $H(E)$ shows that this
set is in bijective correspondence with $\D(E)\cap A,$ via
$\lambda(A,E)\mapsto A,$ and that for any set $B,$ the equivalence
relation $R(A,B,E)$ on $\D(E)\cap A$ corresponds to the
restriction to $H(E)$ of the congruence generated by the
improper equivalence relations on the sets $H(F)$ $(F\in B\cap\D(E)).$
It follows that given $B$ as in~(\ref{finB}), the
improper congruence on $H$ is generated by the finite set of pairs
$\{(\lambda(F,E),\lambda(F',E))\ |\ E\in B;\ F,F'\in A\cap \D(E)\}.$
Conversely, assuming~(\ref{HfgJ}), we may take a finite
generating set $S$ for the improper congruence
on $H$ and let $B=$\hbox{$\{E\ |\ (\exists\,F,F'\in A)$}
$(\lambda(F,E),\lambda(F',E))\in S\},$ and we see that this $B$
witnesses~(\ref{finB}).\qed
\end{pf}

The above lemma, combined with Theorem~\ref{TEfixed}, yields
necessary and sufficient conditions for a category of the form
$J_\fb{cat}$ to have the property we are interested in (last paragraph
of theorem below).
We can also get from it some necessary conditions for this to be
true of an arbitrary small category (first paragraph).

\begin{thm}\label{Tgath}
Let $\fb E$ be a small category, and $J$ the partially
ordered set whose elements are the equivalence classes of
objects of $\fb E$ under the equivalence relation that
relates $E$ and $F$ if there exist morphisms from $E$ to $F$
and from $F$ to $E$ \r(\!\<cf.\ proof of Proposition~\ref{PinjE}\r).
Let $A$ denote the set of minimal elements of $J,$ and $B$
the set of $\!A\!$-critical elements.
Then {\em necessary} conditions for limits over $\fb E$ to respect
direct limits of $\!\fb E\!\<$-sets are {\rm~(i)}~$A$ is finite,
{\rm~(ii)}~every element of $J$ lies above an element of $A,$
{\rm~(iii)}~$B$ is finite, and
{\rm~(iv)}~$B$ gathers $A$ under every element of $J.$

If $\fb E$ is in fact a category formed from a partially ordered set
by the construction $(~)_\fb{cat}$ \r(\!\<equivalently,
if $\fb E \cong J_{\fb{cat}}),$ then the conjunction of these
four conditions is sufficient as well as necessary.
\end{thm}\begin{pf}
The final assertion is immediate from Proposition~\ref{PinjE},
Theorem~\ref{TEfixed}, and Lemma~\ref{Lgath}.

To get the assertion of the first paragraph, suppose
that limits over $\fb E$ respect direct limits of $\!\fb E\!\<$-sets.
Conditions~(i) and~(ii) follow from Proposition~\ref{PinjE}.
Let us write the set of minimal elements of $J$ more distinctively
as $A^{(J)},$ let $A^{(\fb E)}\subseteq\ObE$ be a set of
representatives of these elements, and let
$H^{(J_\fb{cat})}$ and $H^{(\fb E)}$
denote the $\!J_\fb{cat}\!\<$-set and the $\!\fb E\!\<$-set determined
by these respective sets of objects.
Then by Theorem~\ref{TEfixed} our assumption implies that
the trivial congruence on $H^{(\fb E)}$ is finitely generated, while
by Lemma~\ref{Lgath}, the conclusions~(iii) and~(iv) that we want
to prove are equivalent to saying that the same is true
of the trivial congruence on $H^{(J_\fb{cat})}.$

Now there is an obvious functor $R\,{:}\,\ \fb E\to J_\fb{cat}$
taking each object of $\fb E$ to its equivalence class in $J.$
It is easy to see that the composite functor
$H^{(J_\fb{cat})}{\circ}R\,{:}\:\ \fb E\to J_\fb{cat}\to \fb{Set}$
admits a surjective homomorphism
$H^{(\fb E)}\rightarrow H^{(J_\fb{cat})}{\circ}R;$ hence as
the improper congruence on $H^{(\fb E)}$ is finitely generated, the
same is true of the improper congruence on $H^{(J_\fb{cat})}{\circ}R\<,$
and hence, as $R$ is surjective on objects, of the improper congruence
on $H^{(J_\fb{cat})},$ as required.\qed
\end{pf}

\section{Remarks}\label{misc}

As noted in the introduction, given a directed system of
algebras $(A_i)_I,$ understood to be finitary,
one can construct its direct limit by taking
the direct limit of underlying sets and putting an
appropriate algebra structure on this set, essentially because
direct limits respect finite products
of sets, and an algebra structure is given by maps on such product sets.
On the other hand, direct limits do not in general respect infinite
products; indeed, such a product can be thought of as a limit over
$J_\fb{cat}$ where $J$ is an infinite antichain, and such a $J$ does not
satisfy condition~(i) of Theorem~\ref{Tgath}.
So direct limits of infinitary algebras cannot be constructed as in
the finitary case.
An example is
\begin{exmp}\label{0oo}
A directed system of algebras with one $\!\aleph_0\!$-ary
operation, such that the algebra structure cannot be extended
to the direct limit set in any natural way.
\end{exmp}
\begin{pf*}{Details.}
For each positive real number $a$ let $A_a$ be the closed
interval $[0,a]\subseteq\Rset,$ given with the $\!\aleph_0\!$-ary
supremum operation $(x_0,x_1,\dots)\mapsto
\mathrm{sup}(x_0,x_1,\dots).$
These sets form a directed system under inclusion, but the operation
$\mathrm{sup}$ clearly does not extend in a natural way to their
direct limit, $[0,\infty).$
For instance, one has no natural definition of
$\mathrm{sup}(0,1,2,\dots),$
because the map $\iota: \coLm_{a\in\Rset}
(A_a^{\Nset})\rightarrow
(\coLm_{a\in\Rset} A_a)^{\Nset}$ does
not have $(0,1,2,\dots)$ in its image.
It is not hard to show that no extension of
$\mathrm{sup}$ to $[0,\infty)$ makes this set the direct limit
of the algebras $A_a.$
(The uncountability of $\Rset$ is not necessary to
this example; one may replace $\Rset$ with $\Nset.$
I just felt that the supremum function on
real numbers was the more ``important'' example.)
\end{pf*}

We also noted in the introduction that the results of this paper
are specific to $\!\fb{Set}\!$-valued functors, and fail for
functors with other codomains, e.g., $\fb{Set}^\mathrm{op}.$
For another example, let $\fb{Metr}$ be the category of metric
spaces, with distance-nonincreasing maps as morphisms.
Then one has
\begin{exmp}\label{Metr}
A directed system of $\!Z_2\!$-sets $X_0\to X_1\to\cdots,$
in $\fb{Metr}$ such that the map $\iota: \coLm_i (X_i)^{Z_2}\rightarrow
(\coLm_i X_i)^{Z_2}$ is not surjective.
\end{exmp}
\begin{pf*}{Details.}
For each $i,$ let $X_i$ be the set $\{0,1\},$ with $d(0,1)=1/(i+1),$
and with $Z_2$ acting by switching $0$ and $1,$
and let all connecting morphisms be the identity on underlying sets.
Each of
the sets $X_i^{Z_2}$ is empty, so $\coLm_I\,X^{Z_2}_i = \varnothing.$
However, from the metric space axiom $d(x,y)=0\implies x=y$
one sees that the direct limit of this directed system
is the $\!1\!$-point metric space, on which
$Z_2$ acts trivially; thus, $(\coLm_I\,X_i)^{Z_2}$ is nonempty.
\end{pf*}

A type of question related to that considered in
this note arises in sheaf theory.
A sheaf of sets on a topological space $V$ is a certain sort of functor
$(o(V)^\mathrm{op})_\fb{cat}\to\fb{Set},$ where $o(V)$ is the set
of open subsets of $V,$ partially ordered by inclusion;
and an analog of the question we have considered is, ``When does
the global-sections functor commute with direct limits of sheaves?''
But that problem is not actually a case of the problem considered
above, because of the nontrivial form that the direct limit
construction takes for sheaves.
A class of situations where that commutativity holds is obtained
in \cite[Exercise~II.1.11]{RH}, and in greater
generality in \mbox{\cite[Proposition 3.6.3]{AG}}.

\end{document}